\documentclass{amsart}
\usepackage[alphabetic]{amsrefs}
\usepackage[letterpaper, portrait, margin=2.5cm]{geometry}

\usepackage{amsmath,amssymb,amsthm}
\usepackage{mathtools}
\usepackage{tikz-cd}
\usepackage{comment}

\usepackage[table]{xcolor}
\usepackage{colortbl}
\definecolor{c0}{HTML}{000000}
\definecolor{c1}{HTML}{9AB8FD}
\definecolor{c2}{HTML}{E361A0}
\definecolor{c3}{HTML}{FD8F60}
\definecolor{c4}{HTML}{998BF2}
\definecolor{c5}{HTML}{1CCA93}
\definecolor{c6}{HTML}{FEC161}
\definecolor{c7}{HTML}{D40000}

\usepackage{hyperref}

\usepackage{float}
\restylefloat{table}

\usepackage[textwidth=20mm]{todonotes}
\newcommand{\pablo}[2][]{\ifthenelse{\equal{#1}{inline}}{\todo[inline,color=orange]{#2}}{\todo[color=orange]{#2}}}
\newcommand{\akash}[2][]{\ifthenelse{\equal{#1}{inline}}{\todo[inline,color=pink]{#2}}{\todo[color=pink]{#2}}}
\newcommand{\claudio}[2][]{\ifthenelse{\equal{#1}{inline}}{\todo[inline,color=yellow]{#2}}{\todo[color=yellow]{#2}}}
\newcommand{\ella}[2][]{\ifthenelse{\equal{#1}{inline}}{\todo[inline,color=cyan]{#2}}{\todo[color=cyan]{#2}}}
\newcommand{\carmen}[2][]{\ifthenelse{\equal{#1}{inline}}{\todo[inline,color=lime]{#2}}{\todo[color=lime]{#2}}}

\newcommand{\Cb}{\mathbb C}

\newcommand{\Fb}{\mathbb F}

\newcommand{\Pb}{\mathbb P}

\newcommand{\Vb}{\mathbb V}

\newcommand{\Zb}{\mathbb Z}

\DeclareMathOperator{\Spec}{Spec}

\DeclareMathOperator{\GL}{GL}
\DeclareMathOperator{\PGL}{PGL}
\DeclareMathOperator{\SL}{SL}
\DeclareMathOperator{\PSL}{PSL}
\DeclareMathOperator{\SU}{SU}
\DeclareMathOperator{\PSU}{PSU}
\DeclareMathOperator{\Tr}{Tr}

\newcommand{\onto}{\twoheadrightarrow}

\DeclareMathOperator{\Sym}{Sym}

\DeclareMathOperator{\RD}{RD}

\DeclarePairedDelimiter\floor{\lfloor}{\rfloor}

\newtheorem{theorem}{Theorem}[section]

\theoremstyle{definition}
\newtheorem{definition}[theorem]{Definition}
\newtheorem{problem}{Problem}

\newtheorem{example}[theorem]{Example}

\numberwithin{equation}{section}

\title{On the resolvent degree of $\PSU(3,q)$}

\author[Christofferson, Ganguly, G\'{o}mez-Gonz\'{a}les, Kuriyama, and Li]{Pablo Nicolas Christofferson, Akash Ganguly, Claudio G\'{o}mez-Gonz\'{a}les, Ella Kuriyama, \\
and Yihan Carmen Li \\
\ \\
\ with an appendix co-authored by Nawal Baydoun}

\email{christoffersonp@carleton.edu}
\address{Department of Mathematics \& Statistics, Carleton College}

\email{aganguly@umn.edu}
\address{Department of Mathematics, University of Minnesota}

\email{cgonzales@carleton.edu}
\address{Department of Mathematics \& Statistics, Carleton College}

\email{kuriyamae@carleton.edu}
\address{Department of Mathematics \& Statistics, Carleton College}

\email{lic6@carleton.edu}
\address{Department of Mathematics \& Statistics, Carleton College}

\thanks{The authors were supported in part by NSF Grant DMS-2418943.}

\begin{document}
\begin{abstract}
    Resolvent degree ($\RD$) is an invariant of finite groups in terms of the complexity of their algebraic actions. We address the problem of bounding $\RD(G)$ for all finite simple groups using the methods established in \cite{GGSW2024} in terms of $\RD^{\leq d}_{\Cb}$-versality and special points. We give upper bounds on $\RD(\PSU(3,q))$ and $\RD(\PSU(2, q))$ in terms of classical invariant theory.  In the $\PSU(3,q)$ case, stability of low-degree invariants permit an asymptotic bound on $\RD$ growing in $q$.
\end{abstract}
\maketitle

\section{Introduction}

Resolvent degree is a measure of complexity for finite groups motivated by the problem of solving polynomial equations in the simplest manner possible---beyond the usual framework of solutions in radicals. The quadratic, cubic, and quartic formulas, together with other elementary results in Galois theory, show that cyclic groups, $S_3$, $A_4$, and $S_4$ all have $\RD = 1$. Indeed, lesser known solutions to algebraic problems such as Bring's \cite{Bring1786} reduction of the quintic mean that $A_5$ and $S_5$ also have $\RD = 1$. 

The icosahedral solution to the quintic \cite{Klein1884}, together with accompanying notions of versality and accessory irrationality, was paradigmatic of Klein's program to express solutions to algebraic functions in terms of actions of the monodromy group on some variety of minimal dimension. Such a perspective gives rise to the modern notion of $\RD(G)$ (see \cite{FW2019, Reichstein2025}), where $\RD(n) \coloneqq \RD(S_n) = \RD(A_n)$ captures the classical problem of eliminating coefficients from the generic degree $n$ polynomial that has received modern attention as in \cites{Brauer1975,AS1976,Wolfson2020,Sutherland2021,HeberleSutherland2023}. 

While there are currently no effective techniques for establishing lower bounds on $\RD(G)$, the problem of finding upper bounds has received both classical and modern interest for several families of finite groups: $A_n$ and $S_n$; simple factors of the Weyl groups of type $E_6$, $E_7$, and $E_8$; $\PSL(2,7)$ and $\PSL(2,11)$; and the sporadic groups. See \cite{Sutherland2023} for a survey of progress to date. 
Having established bounds on $\RD(G)$ in terms of Jordan--H\"older decompositions, Farb--Wolfson posed:

\begin{problem}[\cite{FW2019}]
    Understand $\RD(G)$ for all finite simple groups $G$.
\end{problem}

The classification of finite simple groups organizes all such $G$ into the cyclic groups of prime order (all have $\RD=1$), the alternating groups $A_n$ (the implicit objects of study in $\RD$ over several centuries), 16 infinite families of Lie type, and the 26 sporadic groups (with bounds on their $\RD$ provided in \cite{GGSW2024}). In this paper, we utilize the method introduced by G\'omez-Gonz\'ales--Sutherland--Wolfson to bound the resolvent degree of subfamilies of the finite unitary groups. Our main result is:
\begin{theorem}[Bounds on the resolvent degree of $\PSU(3,q)$]\label{thm:bounds_of_rd_psu3} We have
\begin{align*}
    \RD(\PSU(3,3)) &\leq 4, & \RD(\PSU(3,4)) &\leq 10, & \RD(\PSU(3,5)) &\leq 17, \\
    \RD(\PSU(3,7)) &\leq 39, & \RD(\PSU(3,8)) &\leq 53, & \RD(\PSU(3,9)) &\leq 69, \\
    \RD(\PSU(3,11)) &\leq 106, & \RD(\PSU(3,13)) &\leq 152, & \RD(\PSU(3,16)) &\leq 236, \\
    \RD(\PSU(3,17)) & \leq 267, & \text{ and } & & \RD(\PSU(3,19)) & \leq 338. \\
\end{align*}
Moreover, for all prime powers $q \geq 23$, we have $\RD(\PSU(3,q))\leq q^2-q-\log_4(q^2-q+6)$.

\end{theorem}

Except for $q=2$, this estimate is a substantial improvement over the best estimates for $\RD(q^3+1)$ in the literature, i.e., the implicit bound on $\RD(\PSU(3,q))$ via the minimal permutation representation (except in the $q=5$ case). We note that, other than the extensive body of work establishing asymptotic bounds on $\RD(n)$, the uniformity of these bounds and process for obtaining them across an entire infinite family of finite non-abelian simple groups is distinctive in the literature, and we hope that this paper and possible sequels might contribute to more systematic strategies for obtaining bounds applying to entire families of Lie type. 

Our main Theorem follows from the construction of a $\PSU(3,q)$-variety of controllable versality, enabled by the well-behaved structure of invariants for a specific representation of $\PSU(3,q)$ serving as the ambient projective space (Theorem \ref{thm:small_invariants}). For the sake of accessibility and self-consistency, we also state a version of \cite[Theorem~4.7]{GGSW2024} (appearing as Theorem \ref{thm:the_game}) which directly gives a criterion for constructing such varieties in the context of all finite simple groups. We also obtain bounds on $\RD$ of $\PSU(2,q) \cong \PSL(2,q)$ for all prime powers $q \leq 125$, using the computer algebra systems SageMath \cite{Sage2025} and GAP \cite{GAP2024}. Several of these calculations reproduce independent results by Nawal Baydoun, and are collected in jointly-written Appendix \ref{appendix:bounds_on_psu2} together with worked examples to describe the method in detail. Finally, while the calculations needed to conclude Theorem \ref{thm:small_invariants} are too extensive to include in full detail, the interested reader is referred to Appendices \ref{appendix:bounds_on_psu3} and \ref{appendix:powers_in_psu3}; the former includes numerics specific to calculations with $q \leq 125$ and the latter contains the power map data needed to recover said results.

\subsection{Acknowledgments}
The authors thank Alexander Sutherland and Jesse Wolfson for many useful comments on an early draft. The authors are especially appreciative to Nawal Baydoun for sharing unpublished notes and preliminary results for $G = \PSL(2,q)$ which clarified the process underlying Theorem \ref{thm:the_game}; in addition, she and Jonathan Delgado were invaluable in their work independently verifying calculations.
\section{Background}

\subsection{Representation theory}

In this section we recall the requisite notions from representation and character theory; while many standard sources exist (e.g., \cite{Isaacs1994}), we aim for consistency of language and relative self-contained-ness in describing the associated algorithm. All representations are of finite groups on finite-dimensional $\Cb$ vector spaces.

For a finite group $G$ with representation $W$, we are interested in the representation $\Sym^k W^\ast$, i.e., the space of degree $k$ homogeneous polynomials on $W$, and its invariants. The \emph{Molien series} of $W$ is
\[ M_{G,W}(t) = \sum_{k \geq 0} m_k(G,W) \, t^k, \qquad \text{ where } \qquad m_k(G,W) = \dim \left( \Sym^k W^\ast \right)^G, \]
where we simply write $M_W(t)$ and $m_k(W)$ when $G$ is clear from context. In other words, the $k$th coefficient of $M_W(t)$ measures the space of degree $k$ homogeneous $G$-invariant polynomials on $W$. \emph{Molien's formula} provides a closed form expression for this series: \[ M_{G,W}(t) = \frac{1}{|G|} \sum_{g \in G} \det\left(1 - tg|_{W^\ast}\right)^{-1}. \]
While various computer algebra systems (e.g., \cite{GAP2024}) are capable of determining $M_W(t)$ for representations of sufficiently small groups, for our purposes, the entire Molien series is unnecessary and prohibitively costly to compute; a more effective approach is via character theory. 

\bigskip
The character $\chi_W$ of a representation $\rho: G \to \GL(W)$ is given by $\chi_W(g) = \Tr(\rho(g))$. This $\Cb$-valued function, which is constant on conjugacy classes of $G$, classifies $\rho$ up to isomorphism. The formula
\begin{equation}\label{eqn:invariant_dimension}
\dim W^G = \frac{1}{|G|} \sum_{g \in G} \chi_W(g) = \sum_{i=1}^r \frac{\chi_W(g_i)}{|C(g_i)|},
\end{equation}
where we enumerate class representatives $g_1,\dots,g_r$ and write $C(g_i)$ for the corresponding centralizer, with
\begin{equation}\label{eqn:character_identities}
\chi_{W^\ast}(g) = \overline{\chi_W(g)} \quad \text{ and } \quad \chi_{\Sym^k W}(g) = \frac{1}{k} \sum_{i=1}^k \chi_{\Sym^{k-i}(W)}(g) \chi_W(g^i) 
\end{equation}
together allow us to compute the terms $m_k(G,W)$ directly for small values of $k$.

A representation of $G$ is called \emph{irreducible} if it has no proper non-trivial subrepresentations; a fundamental result in representation theory is that a finite group has exactly as many non-isomorphic irreducible representations as it does conjugacy classes. This fact leads naturally to the notion of a \emph{character table}, which lists all the irreducible representations of $G$ on one axis and conjugacy classes representatives of $G$ on the other: each entry gives the value of a particular irreducible character evaluated at a specific conjugacy class. Many computer algebra systems (e.g., \cite{GAP2024}) contain (or are capable of computing) character tables for sufficiently small groups. For our purposes, considering the smallest non-trivial irreducible representations (i.e., not the entire table) is enough to give the best bounds on $\RD(G)$ attainable by our approach. 

\bigskip
The process of bounding $\RD(G)$ described in \cite{GGSW2024} and employed here makes use of \emph{projective} representations, i.e., homomorphisms $G \to \PGL(W)$. Any (linear) representation of $G$ passes to such an action via the quotient map $\GL(W) \onto \PGL(W)$; however, these are not the only projective representations of $G$. Indeed, there always exists a cover $\Gamma \onto G$, called a \emph{Schur cover of $G$}, such that every projective representation of $G$ lifts to a linear representation of $\Gamma$, i.e., the following diagram commutes:
\[
\begin{tikzcd}
    \Gamma \arrow[r, dashed] \arrow[d, two heads] & \GL(V) \arrow[d, two heads] \\ 
    G \arrow[r] & \PGL(W).
\end{tikzcd}
\]
Note that a Schur cover is a central extension of $G$ by its \emph{Schur multiplier} $M(G) \coloneqq H_2(G;\Zb)$. When $G$ is perfect (in particular, for $G$ a finite simple group), $\Gamma$ is unique up to isomorphism and referred to as \emph{the} Schur cover of $G$; we shall do so henceforth. For a more thorough discussion, see \cite[Chapter~11]{Isaacs1994}. 

Identifying $G$-invariant degree $k$ homogeneous polynomials on $\Pb(W)$ is tantamount to computing
\[ \dim \Sym^k_\Gamma(W^\ast). \]
The zero sets of such polynomials define $G$-invariant projective hypersurfaces in $\Pb(W)$. 

\subsection{Resolvent degree of finite groups}\label{subsec:RD_G}

In what follows, all varieties are over $\Cb$. Given our objectives to bound the $\RD$ of finite simple groups, by \cite[Theorem~1.2~and~1.3]{Reichstein2025}, there is no loss of generality.

The resolvent degree of a branched cover $\pi: X \dashrightarrow Y$ (a generically finite, dominant rational map) is defined in \cite{FW2019} in terms of \emph{formulas} for $\pi$, i.e., finite towers of branched covers
\[ Y_r \dashrightarrow Y_{r-1} \dashrightarrow \cdots \dashrightarrow Y_1 \dashrightarrow Y_0 = Y, \]
such that the composite $Y_r \dashrightarrow Y$ factors through a branched cover $Y_r \dashrightarrow X$. For a given formula, we are concerned with the largest \emph{essential dimension} of each $Y_i \dashrightarrow Y_{i-1}$; the resolvent degree $\RD(\pi)$ is the minimal such uniform bound over all formulas for $\pi$. This notion precisely quantifies the spirit of Hilbert's 13th problem (and sextic, octic, and nontic conjectures), where the essential dimension of the intermediate covers are understood as the number of variables necessary to write down a formula solving the associated algebraic problems. In particular, because every branched cover is a formula for itself (a tower of length 1), we always have that $\RD(\pi) \leq \dim X$.

\begin{definition}\label{def:rd_group} Let $G$ be a finite group. Then
\[ \RD(G) \coloneqq \sup\{ \RD(X \dashrightarrow X/G): X \text{ a primitive, faithful $G$-variety} \}. \]
\end{definition}
So $\RD(G)$ measures the complexity of algebraic actions of $G$; in particular, this construction is designed so that $\RD(S_n)$ is the minimal number of variables needed to solve the generic degree $n$ polynomial. As is tradition, we write $\RD(n) \coloneqq \RD(S_n)$. There is extensive literature dating back to at least \cite{Tschirnhaus1683} in obtaining upper bounds on $\RD(n)$; see \cite[Appendix~B]{Wolfson2020} for a historical overview.

Farb--Wolfson show that, for any short exact sequence $1 \to N \to G \to H \to 1$ of finite groups, we have
\[ \RD(G) \leq \max\{\RD(N),\RD(H)\}. \]
More generally, the resolvent degree of a finite group is bound by the resolvent degree of simple factors in its Jordan--H\"older decomposition. The natural problem of understanding $\RD(G)$ for any finite simple group $G$ is that which animates this paper. For more details and a recent summary of known bounds, see \cite{Sutherland2023}. 

\bigskip Note that, in what follows, we write $\mu(G)$ for the degree of the minimal permutation representation of $G$, i.e., the cardinality of the smallest sets on which $G$ can act non-trivially. We make use of a weakened version of \cite[Theorem~4.7]{GGSW2024}:

\begin{theorem}[Bounds on resolvent degree via invariant theory]\label{thm:the_game} Let $G$ be a finite non-abelian simple group with a non-trivial irreducible projective representation $\Pb(V)$. Let $f_1,\dots,f_r$ be $G$-invariant homogeneous polynomials on $V$ with respective degrees $d_1 \leq \dots \leq d_r$, such that the $f_i$ are algebraically independent from the lower degree invariants. If $d_1 \cdots d_r < \mu(G)$ and $\RD(d_1 \cdots d_r) \leq \dim V - 1 - r$, then
\[ \RD(G) \leq \dim \Vb(f_1,\dots,f_r) \leq \dim V - 1 - r. \]
\end{theorem}

Essential to the theory of resolvent degree is that of \emph{versality} and its generalizations, which we describe functionally but avoid defining for the sake of brevity (see \cite[Section 3]{GGSW2024} for a thorough accounting). In particular, versality is a sufficient condition to witness the supremum of Definition \ref{def:rd_group}: by \cite[Proposition 3.10]{FW2019}, if $X$ is a versal $G$-variety then $\RD(G) = \RD(X \dashrightarrow X/G)$. On the other hand, as was known classically, versality is not necessary; Klein showed that $\RD(A_5) = 1$ using the projective representation $A_5 \curvearrowright \Pb_\Cb^1$ while also exhibiting its non-versality \cite[Chapter 5, \S11]{Klein1884}. Such phenomena gives rise to the idea of \emph{versality up to some accessory irrationality}, with the additional requisite square root in Klein's solution of the quintic as the seminal example. 

These broadened notions of versality expand the class of $G$-varieties available for analysis in understanding $\RD(G)$. Various classes of accessory irrationalities appear prominently in work on resolvent problems, where the modern perspective is given formally in \cite{FKW2023} using the language of branched covers. In \cite{GGSW2024}, the authors (using the dual language of field extensions) establish the technique used herein to bound $\RD(G)$ via constructing explicit $G$-varieties that are \emph{$\RD^{\leq d}_{\Cb}$-versal}, corresponding to a class of accessory irrationality permitting auxiliary solutions to equations with resolvent degree bounded by $d$. 

The relevance of this class, first studied in \cite{AS1976}, to our present endeavors follows from \cite[Proposition 3.3]{GGSW2024}:
\begin{equation}\label{eqn:rd_g_bound} \RD(G) = \min_{d \geq 0} \left\{ \max\{d, \dim X\} \, \big| \, X \text{ is an $\RD^{\leq d}_{\Cb}$-versal $G$-variety} \right\}. \end{equation}

\begin{proof}[Sketch of Theorem \ref{thm:the_game}]
We mirror the argument of \cite[Theorem~4.7]{GGSW2024}--- which we urge the reader to compare with what appears here for a thorough accounting of versality, $G$-torsors, and twists---the core of which applies for any non-abelian finite simple group. Our goal is to illustrate the constraints giving rise to the algorithm demonstrated in Appendix \ref{appendix:bounds_on_psu2}. 
Setting $d = \RD(d_1 \cdots d_r)$, the proof amounts to showing that 
\[ X = \Vb(f_1,\dots,f_r) \subseteq \Pb(V) \]
is $\RD^{\leq d}_{\Cb}$-versal and applying \eqref{eqn:rd_g_bound}. This is accomplished through the criterion of \cite[Theorem~3.9]{GGSW2024}, which includes verifying generic freeness of the action $G \curvearrowright X$ and a density condition in terms of twists $\prescript{T}{}{X}$ for every $G$-torsor $T \to \Spec K$ with $K$ finitely generated over $\Cb$. 

In light of \cite[Lemma~2.15]{GGSW2024}, which demonstrates that irreducibility is a sufficient condition for generic freeness of the $G$-action, the first criterion is straightforward: since $X$ has at most $d_1 \cdots d_r$ many components and the action of $G$ permutes them, we conclude that $X$ is irreducible since $d_1 \cdots d_r < \mu(G)$. The second criterion is showing that $K^{(d)}$-points are dense in every twisted form of $X$, where $K^{(d)}$ is defined with respect to some algebraic closure as a compositum of intermediate field extensions with resolvent degree bounded by $d$ as in \cite[Definition~2.17]{GGSW2024}. While the details are outside the scope of this paper, the density argument follows nearly verbatim by that of \cite[Theorem~4.7]{GGSW2024} (specifically the case where, in their notation, $Y_G = \Pb(V_G)$) via \cite[Lemma~14.4]{Reichstein2025} and the assumption $d \leq \dim V - 1 - r$.
\end{proof}

We include several concrete examples of applying Theorem \ref{thm:the_game}, phrased as an algorithm which allows us to methodically cut down $\RD(G)$ using invariants of a projective representation of $G$, in Appendix \ref{appendix:bounds_on_psu2}. From this perspective, the codimension of the intersection---the amount by which we are able to reduce the bound on $\RD(G)$---is constrained by the problem of finding $K^{(d)}$-dense points on a variety of high degree.

\clearpage
\section{Resolvent degree of \texorpdfstring{$\PSU(3,q)$}{PSU3q}}

\subsection{Invariant forms}\label{subsec:invariant_forms}
Throughout this paper, we write $\PSU(n,q)$ where some use the notation $\PSU(n,q^2)$. The conjugacy classes of $\PSU(3,q)$ are given in \cite{FS1973} in terms of Jordan normal forms (with corrections mentioned in \cite{Orevkov2013}). These classes (indeed, the character tables) admit a highly uniform description, organized into types according to centralizer sizes, where indices ($k, \ell,$ and $m$) are used to describe classes of a similar type and some of their relations. Writing $\delta_{i,j}$ for the Kronecker delta, we have

\begin{center}
\footnotesize
\renewcommand{\arraystretch}{1.75}
\begin{tabular}{r|c|c|c|c|c|c|c|c|c|}
Type & $C_1$ & $C_2$ & $C_3^\ell$ & $C_4^k$ & $C_5^k$ & $C_6'$ & $C_6^{k,\ell,m}$ & $C_7^k$ & $C_8^k$ \\
\hline
$|\text{Centralizer}|$ & $|G|$ & $\tfrac{q^3(q+1)}{d}$ & $q^2$ & $\tfrac{(q-1)q(q+1)^2}{d}$ & $\tfrac{q(q+1)}{d}$ & $(q+1)^2$ & $\tfrac{(q+1)^2}{d}$ & $\tfrac{q^2-1}{d}$ & $\tfrac{q^2-q+1}{d}$ \\
\hline
\# of Classes & $1$ & $1$ & $d$ & $\tfrac{q+1}{d}-1$ & $\tfrac{q+1}{d}-1$ & $\delta_{3,d}$ & $\tfrac{q^2-q+1-d}{6d}$ & $\tfrac{q^2-q+1-d}{2d} - \delta_{1,d}$ & $\tfrac{q^2-q+1-d}{3d}$ 
\end{tabular}
\end{center}

\noindent Note that $d = \gcd(3,q+1)$ is the Schur multiplier of $\PSU(3,q)$, i.e., $\PSU(3,q) = \SU(3,q)$ if $q \not \equiv 2 \bmod 3$ and 
\[ \PSU(3,q) = \SU(3,q) \Big/ \left\langle \small \left(\begin{smallmatrix}
    \omega & 0 & 0 \\
    0 & \omega & 0 \\
    0 & 0 & \omega \\
\end{smallmatrix}\right) \right\rangle \]
otherwise, where $\omega \in \Fb_{q^2} \setminus \Fb_{q}$ is a primitive third root of unity. Indeed, the conjugacy class $C'_6$ is nontrivial if and only if $q \equiv 2 \bmod 3$. A series of routine calculations in terms of conjugacy class representatives yield the power tables displayed in Appendix \ref{appendix:powers_in_psu3}.

\bigskip
While any representation of the Schur cover $\SU(3,q)$ will satisfy the hypothesis of Theorem \ref{thm:the_game}, we are interested in its smallest nontrivial representations in order to obtain the best bounds on resolvent degree. While we used Sage \cite{Sage2025} and GAP \cite{GAP2024} to study all nontrivial irreducible representations for small values of $q$, our analysis here focuses on the smallest nontrivial irreducible representation of $\PSU(3,q)$, which lifts to a representation of $\SU(3,q)$ and always produced the best possible estimates in our calculations. The character of this representation is denoted by $\chi_{qs}$ in \cite{FS1973} and has degree $q^2-q$. While we also considered a family of representations with degree $q^2-q+1$, we verified that $\chi_{qs}$ produced a better (or equal) bound on resolvent degree for all $q \leq 197$ and, given the experimental data and concrete results we achieved in terms of $\chi_{qs}$, we chose to restrict our attention to the smaller representation. Note that all other irreducible characters have degree growing cubically in $q$; even with plentiful invariants in degree $2$, such representations cannot hope to outdo the bounds furnished using the methods of this paper.

Henceforth, we write $V$ for the unique irreducible representation of $\PSU(3,q)$ of degree $q^2-q$, which has the character
\begin{center}
\renewcommand{\arraystretch}{1.75}
\begin{tabular}{r|c|c|c|c|c|c|c|c|c}
Type & $C_1$ & $C_2$ & $C_3^\ell$ & $C_4^k$ & $C_5^k$ & $C_6'$ & $C_6^{k,\ell,m}$ & $C_7^k$ & $C_8^k$  \\
\hline
$\chi_V$ & $q(q-1)$ & $-q$ & $0$ & $1-q$ & $1$ & $2$ & $2$ & $0$ & $-1$. \\
\end{tabular}
\end{center}
This representation is particularly amenable to our analysis because it is self-dual ($V \cong V^\ast$, i.e., $\chi_V = \overline{\chi_V}$) and constant on conjugacy classes of the same type.

\bigskip
Our main Theorem in this section concerns the representations $\Sym^2 V^\ast$, $\Sym^3 V^\ast$, and $\Sym^4 V^\ast$, whose characters are too complicated to write here in full generality, particularly because of their dependence of the value of $q$ modulo $72$. We include in Appendix \ref{appendix:powers_in_psu3} the power maps $g \mapsto g^k$ of $\PSU(3,q)$ for $2 \leq k \leq 4$, from which the desired invariants of these representation can be recovered using \eqref{eqn:character_identities}. All told, we have:

\begin{theorem}
\label{thm:small_invariants}
Let $q$ be a prime power, and let $V$ be the unique irreducible representation of $\PSU(3,q)$ of degree $q^2-q$. Then we have $M_V(t) = 1 + m_4(V) t^4 + O(t^5),$ where
\[ m_4(V) = \left\{ \begin{array}{ll}
    \frac{1}{6}(q-1) & \text{if } q \equiv 1 \bmod 6 \\[0.4em]
\frac{1}{6}(q+10) & \text{if } q \equiv 2 \bmod 6 \\[0.4em]
\frac{1}{6}(q-3) & \text{if } q \equiv 3 \bmod 6 \\[0.4em] 
\frac{1}{6}(q+2) & \text{if } q \equiv 4 \bmod 6 \\[0.4em]
\frac{1}{6}(q+7) & \text{if }q \equiv 5 \bmod 6.
\end{array} \right. \]
\end{theorem}

\begin{proof}
Using Equations \eqref{eqn:character_identities} to calculate the characters $\chi_{\Sym^k V}$ for $2 \leq k \leq 4$, we compute
\[ \dim (\Sym^k V^\ast )^{\PSU(3,q)} = \frac{1}{|\PSU(3,q)|} \sum_{g \in \PSU(3,q)} \chi_{\Sym^k V^\ast}(g), \]
via Equation \eqref{eqn:invariant_dimension}. Using that $\chi_V$ is not only a class function but constant on the \emph{types} of conjugacy classes described above---together with the given data regarding the orders of centralizers, number of conjugacy classes of each type, and power tables in Appendix \ref{appendix:powers_in_psu3}---we find that
\[ \dim (\Sym^2 V^\ast)^{\PSU(3,q)} = 0 = \dim (\Sym^3 V^\ast)^{\PSU(3,q)} \]
holds generally, as well as the stated formula for $m_4(V) = \dim (\Sym^4 V^\ast)^{\PSU(3,q)}$.

While these calculations depend on $q \bmod 36$ and are prohibitively cumbersome to express, we include such a computation to demonstrate the process; others follow similarly. Firstly, note that for any finite $G$ with conjugacy classes $g_1,\dots,g_r$ and irreducible representation $W$, we have
\[ 
\dim (\Sym^2 W)^G = \sum_{i=1}^r \frac{\chi_{\Sym^2 W}(g_i)}{|C(g_i)|} = \frac{1}{2} \sum_{i=1}^r \left( \frac{\chi_{W}(g_i)^2}{|C(g_i)|} + \frac{\chi_{W}(g_i^2)}{|C(g_i)|} \right).
\]
If we suppose further that $W$ is self-dual, then we can apply Schur's Lemma to conclude the simpler formula:
\[ 
\dim (\Sym^2 W)^G = \frac{1}{2} + \frac{1}{2} \sum_{i=1}^r \frac{\chi_W(g_i^2)}{|C(g_i)|}.
\]

Returning to the proof, consider the case when $q \equiv 1,3 \bmod 6$ (in particular, $d=1$). Grouping classes by type and employing the power tables from Appendix \ref{appendix:powers_in_psu3}, we compute:
\[ \begin{split}
\dim (\Sym^2 V^\ast)^{\PSU(3,q)} & = \frac{1}{2} + \frac{1}{2} \Bigg( 
    \underbrace{\frac{\chi_V(C_1)}{q^3(q+1)^2(q-1)(q^2-q+1)} }_{C_1} +
    \underbrace{\frac{\chi_V(C_2)}{q^3(q+1)} }_{C_2} + 
    \underbrace{\frac{\chi_V(C_3^\ell)}{q^2} }_{C_3^\ell} \\
    & \quad + \underbrace{\frac{\chi_V(C_1)}{(q-1)q(q+1)^2}+(q-1)\frac{\chi_V(C_4^k)}{(q-1)q(q+1)^2} }_{C_4^k} + 
    \underbrace{\frac{\chi_V(C_2)}{q(q+1)}+(q-1)\frac{\chi_V(C_5^k)}{q(q+1)} }_{C_5^k} \\
    & \quad + \underbrace{\frac{q-1}{2} \frac{\chi_V(C_4^k)}{(q+1)^2} + \frac{(q-1)(q-3)}{6} \frac{\chi_V(C_6^{k,\ell,m})}{(q+1)^2} }_{C_6^{k,\ell,m}} +
    \underbrace{\frac{q+1}{2} \frac{\chi_V(C_4^k)}{q^2-1} + \frac{(q+1)(q-3)}{2} \frac{\chi_V(C_7^k)}{q^2-1} }_{C_7^k} \\
    & \quad + \underbrace{\frac{q(q-1)}{3} \frac{\chi_V(C_8^k)}{q^2-q+1} }_{C_8^k} 
\Bigg) \\
& = \frac{1}{2} + \frac{1}{2} \Bigg( 
    \frac{q(q-1)}{q^3(q+1)^2(q-1)(q^2-q+1)} + \frac{-q}{q^3(q+1)} + \frac{q(q-1)}{(q-1)q(q+1)^2} + \\
    & \quad + (q-1)\frac{1-q}{(q-1)q(q+1)^2} + \frac{-q}{q(q+1)} + (q-1)\frac{1}{q(q+1)} + \frac{q-1}{2} \frac{1-q}{(q+1)^2} \\
    & \quad + \frac{(q-1)(q-3)}{6} \frac{2}{(q+1)^2} + \frac{q+1}{2} \frac{1-q}{q^2-1} + \frac{q(q-1)}{3} \frac{-1}{q^2-q+1}
\Bigg) \\
& = \frac{1}{2} + \frac{1}{2} (-1) = 0.
\end{split} \]

Given the general formulas for characters of symmetric powers,
\begin{gather*}
    \chi_{\Sym^3 W}(g) = \frac{1}{6}\left(\chi_{W}(g)^3+3\chi_{W}(g)\chi_{W}(g^2) + 2\chi_{W}(g^3)\right) \qquad \text{ and } \\
    \chi_{\Sym^4 W}(g) = \frac{1}{24} \left(\chi_{W}(g)^4+6\chi_{W}(g)^2\chi_{W}(g^2) + 3\chi_{W}(g^2)^2 + 8\chi_{W}(g)\chi_{W}(g^3) + 6\chi_{W}(g^4) \right), \\
\end{gather*}
further calculations can be accomplished by a similar strategy, where sums are carried out one at a time for each summand in such expansions of $\chi_{\Sym^k V^\ast}(g)$. To evaluate the sum for a given term, since $k \leq 4$, we need only compare products of something depending on the type of $g$ (the order of its centralizer and possibly $\chi_V(g)$) with $\chi_V(g^k)$ and so each sum can be evaluated using only $\chi_V$ and the data in Appendix \ref{appendix:powers_in_psu3}.
\end{proof}

\subsection{Bounding \texorpdfstring{$\RD(\PSU(3,q))$}{RD(PSU(3,q))}}

We make use of the minimal permutation degree \cite[Table~4]{GMPS2015}:
\begin{equation}\label{eqn:minimal_perm_psu3}
\mu(\PSU(3,q)) = \left\{ \begin{array}{ll}
    q^3 + 1 & q \not = 5 \\
    50 & q = 5.
\end{array} \right.
\end{equation}

We are now equipped to conclude our main result.

\begin{proof}[Proof of Theorem \ref{thm:bounds_of_rd_psu3}]
The cases for $q \leq 19$ are given explicitly in Table \ref{table:PSU3q_bounds}. We will take $q \geq 23$, where there is always a sufficient quantity of algebraically independent quartic hypersurfaces to complete the argument in light of Theorem \ref{thm:small_invariants}. We choose the largest $r$ such that $4^r+r \leq q^2-q+6$; note that, in the given range, this is either $r = \floor{\log_4(q^2-q+6)}$ or $\floor{\log_4(q^2-q+6)}-1$, and we always have $r \geq 4$.

Consider an intersection $X$ of $r$ quartic hypersurfaces in $\Pb(V)$. The bound established by Sutherland \cite{Sutherland2021} gives $\RD(n) \leq n-7$ for $n \geq 109$, so we always have $\RD(4^r) \leq 4^r-7$. Thus
\[ \RD(4^r) \leq q^2 - q - 1 - r = \dim V - 1 - r,  \]
since we can choose the quartics to be algebraically independent. 
Moreover, since $\mu(\PSU(3,q)) = q^3+1 > q^2-q+6$ for all $q \geq 2$, we have $4^r < \mu(\PSU(3,q))$ and so Theorem \ref{thm:the_game} applies. We conclude that
\[ \RD(\PSU(3, q)) \leq q^2-q - 1 - r. \]
Given our choice of $r$, the desired formula follows.
\end{proof}

Note that Hamilton \cite{Hamilton1836} showed that $n-\RD(n)$ grows without bound as $n \to \infty$, and subsequent efforts have established asymptotic upper bounds on the growth of $\RD(n)$. While one could take more care in the above proof to replace $\RD(n) \leq n-7$ with an asymptotic bound, such as those in \cite{Sutherland2021}, the savings given current estimates would be minor and much more complicated to state.

\newpage
\appendix

\section{Explicit bounds on \texorpdfstring{$\RD(\PSU(2,q))$}{RD(PSU(2,q))}, joint with Nawal Baydoun}\label{appendix:bounds_on_psu2}

The representations of the groups $\PSU(2,q) \cong \PSL(2,q)$, which are simple whenever $q \geq 4$, are of classical interest \cite{Jordan1907, Schur1907}. The (oftentimes trivial) Schur cover of $\PSL(2,q)$ is $\SL(2,q)$, except for when $q=4$ and $q=9$ in which case $\PSL(2,q) = \SL(2,q)$ has Schur multiplier given by $\Zb/2\Zb$ and $\Zb/6\Zb$, respectively \cite[Theorem~7.1.1]{Karpilovsky1987}; these cases are already well-studied, since $\PSL(2,4) = A_5$ and $\PSL(2,9) = A_6$. We note the smallest permutation representation of $G$, as in \cite{GMPS2015}:
\begin{equation}\label{eqn:minimal_perm_psu2}
\mu(\PSU(2,q)) = \left\{ \begin{array}{ll}
    q + 1 & q \not = 5, 7, 9, 11 \\
    5, 7, 6, 11 & q = 5, 7, 9, 11, \text{ respectively.}
\end{array} \right.
\end{equation}
In addition to applying Theorem \ref{thm:the_game}, we also need $\mu(G)$ because \cite[Lemma~3.13]{FW2019} gives the bound
\begin{equation}\label{eqn:bound_by_permutation_rep}
    \RD(G) \leq \RD(\mu(G)).
\end{equation}

This section summarizes the process of applying Theorem \ref{thm:the_game} on an irreducible projective representation of $\PSU(2,q)$ for small values of $q$, exhaustively choosing as many invariants as possible before one (or both) of the premises $d_1 \cdots d_r < \mu(\PSU(2,q))$ or $\RD(d_1 \cdots d_r) \leq \dim V - 1 - r$ are violated. We note that, in the context of resolvent degree, $\PSL(2,7)$ was studied classically and refer the reader to \cite{FKW2023}, which establishes $\RD(\PSL(2,7)) = 1$. Indeed, \cite{Klein1878} showed $\RD(\PSL(2,11))) \leq 3$ and \cite{Sutherland2023} gives the improved bound $\RD(\PSL(2,11)) \leq 2$ using the full power of \cite[Theorem~4.7]{GGSW2024}. 

\begin{example}[\cite{Klein1878}]
$\SL(2,7)$ has an irreducible representation $V$ of degree $3$, passing to a projective representation of $\PSL(2,7)$ of dimension 2, and so $\RD(\PSL(2,7)) \leq 2$. Next, a routine calculation gives
\[ M_{V}(t) =  1 + t^4 + t^6 + t^8 + \cdots \] 
Taking $f_1$ to be a degree $4$ invariant, we apply Theorem \ref{thm:the_game} with the curve $\Vb(f_1)$ and find $\RD(\PSL(2,7)) = 1$.
\end{example}

\begin{example}
$\SL(2,13)$ has an irreducible representation $V$ of degree $6$, passing to a projective representation of dimension 5. We observe that
$M_{V}(t) =  1 + t^4 + 2t^8 + 2t^{10} + \cdots $. Setting $f_1$ to be a degree $4$ invariant, Theorem \ref{thm:the_game} applies to the hypersurface $\Vb(f_1)$ and hence $\RD(\PSL(2,13)) \leq 4$. We have additional invariants in degree $8$; if we choose $f_2$ to be such an invariant outside the span of $(f_1)^2$, then $\Vb(f_1,f_2)$ will have dimension $3$. However, we can no longer guarantee irreducibility of the variety $\Vb(f_1,f_2)$, since
\[ \mu(\PSL(2,13)) = 14 \leq 4 \cdot 8 = 32. \] 
Moreover, as the best known bound on $\RD(32)$ is $26$, we cannot extract meaningful information from \eqref{eqn:rd_g_bound}, i.e., the second premise of Theorem \ref{thm:the_game} does not apply. In either case, the algorithm terminates, and we conclude that $\RD(\PSL(2,13)) \leq 4$. 
\end{example}

\begin{example}
$\SL(2,71)$ has an irreducible representation $V$ of degree $35$, passing to a projective representation of dimension 34. We observe that
$M_{V}(t) =  1 + 3t^4 + 2t^5 + 40t^{6} + \cdots $. Setting $f_1,f_2$ to be two linearly independent degree 4 invariants, Theorem \ref{thm:the_game} applies to $\Vb(f_1,f_2)$ and hence $\RD(\PSL(2,71))\leq 32$. We have an additional independent invariant $f_3$ in degree 4, however given the current best upper bounds on resolvent degree we can not conclude $\RD(64) \leq 31$ and so we cannot further apply Theorem \ref{thm:the_game}, despite irreducibility of $\Vb(f_1,f_2,f_3)$.
\end{example}

In particular, the final example illustrates how, as $q$ grows, the primary obstruction to our method lies in verifying the density of $K^{(d)}$-points for arbitrary twists rather than in verifying irreducibility of $\Vb(f_1,\dots,f_r)$. 

\bigskip
We summarize this process for many small values of $q$ in the following Table. All representations selected are among the smallest nontrivial irreducible representations of the Schur cover, i.e., the smallest nontrivial projective representations of $\PSU(2,q)$; these were found to give the best bound on resolvent degree. Each row includes the bound given by applying Theorem \ref{thm:the_game} as above, the degrees of the invariants used in the application, the cardinality of the minimal permutation representation, and the associated bound on $\RD$. On the coloring:
\begin{itemize}
    \item If the application of Theorem \ref{thm:the_game} does \emph{not} improve on the best bounds in the literature---in most cases these are implicit, as $\dim V-1$ or $\RD(\mu(\PSU(2,q)))$---the row is colored \textbf{purple}.
    \item Otherwise, the row is \textbf{orange} if the obstruction to applying Theorem \ref{thm:the_game} with additional invariants includes irreducibility and \textbf{pink} if the only obstruction comes from $\RD(d_1 \cdots d_r)$ growing too large. 
\end{itemize}

\newpage
\begin{table}[h!]
\centering
\begin{tabular}{|c||c|c|c|c|c|}
\hline
$q$ & $\dim V$   & 
Bound by Thm. \ref{thm:the_game} & Degree of Invariants & $\mu(G)$ & Bound by $\RD(\mu(G))$ \\
\hline
\hline
\rowcolor{c4} 5   & 2          & 1              & None            & 5 & 1 \\
\hline
\rowcolor{c4} 7   & 3          & 1              & 4               & 7 & 3 \\
\hline
\rowcolor{c4} 8   & 7          & 4              & 2, 4            & 9 & 4 \\
\hline
\rowcolor{c4} 11  & 5          & 3              & 3               & 11 & 6 \\
\hline
\rowcolor{c3} 13  & 6          & 4              & 4               & 14 & 9 \\
\hline
\rowcolor{c4} 16  & 15         & 12             & 2, 4            & 17 & 12 \\
\hline
\rowcolor{c3} 17  & 8          & 6              & 8               & 18 & 13 \\
\hline
\rowcolor{c2} 19  & 9          & 7              & 3               & 20 & 15 \\
\hline
\rowcolor{c3} 23  & 11         & 9              & 4         & 24 & 18 \\
\hline
\rowcolor{c3} 25  & 12         & 10             & 4      & 26 & 20 \\
\hline
\rowcolor{c3} 27  & 13         & 10             & 3, 4      & 28 & 22 \\
\hline
\rowcolor{c3} 29  & 14         & 12             & 4      & 30 & 24 \\
\hline
\rowcolor{c3} 31  & 15         & 12             & 4, 4      & 32 & 26 \\
\hline
\rowcolor{c4} 32  & 31         & 27             & 2, 4, 4   & 33 & 27 \\
\hline
\rowcolor{c3} 37  & 18         & 15             & 4, 4      & 38 & 32 \\
\hline
\rowcolor{c2} 41  & 20         & 18             & 4      & 42 & 36 \\
\hline
\rowcolor{c3} 43  & 21         & 18             & 3, 4      & 44 & 38 \\
\hline
\rowcolor{c3} 47  & 23         & 20             & 4, 4      & 48 & 42 \\
\hline
\rowcolor{c3} 49  & 24         & 21             & 4, 4   & 50 & 44 \\
\hline
\rowcolor{c3} 53  & 26         & 23             & 4, 4      & 54 & 48 \\
\hline
\rowcolor{c2} 59  & 29         & 26             & 3, 4      & 60 & 54 \\
\hline
\rowcolor{c3} 61  & 30         & 27             & 4, 4   & 62 & 56 \\
\hline
\rowcolor{c4} 64  & 63         & 59             & 2, 4, 4   & 65 & 59 \\
\hline
\rowcolor{c2} 67  & 33         & 30             & 3, 4      & 68 & 62 \\
\hline
\rowcolor{c2} 71  & 35         & 32             & 4, 4      & 72 & 66 \\
\hline
\rowcolor{c2} 73  & 36         & 33             & 4, 4   & 74 & 68 \\
\hline
\rowcolor{c2} 79  & 39         & 36             & 4, 4      & 80 & 74 \\
\hline
\rowcolor{c2} 81  & 40         & 37             & 4, 4   & 82 & 76 \\
\hline
\rowcolor{c2} 83  & 41         & 38             & 3, 4      & 84 & 78 \\
\hline
\rowcolor{c2} 89  & 44         & 41             & 4, 4   & 90 & 84 \\
\hline
\rowcolor{c2} 97  & 48         & 45             & 4, 4   & 98 & 92 \\
\hline
\rowcolor{c2} 101 & 50         & 47             & 4, 4   & 102 & 96 \\
\hline
\rowcolor{c2} 103 & 51         & 48             & 4, 4      & 104 & 98 \\
\hline
\rowcolor{c3} 107 & 53         & 49             & 3, 4, 4   & 108 & 102 \\
\hline
\rowcolor{c2} 109 & 54         & 51             & 4, 4   & 110 & 103 \\
\hline
\rowcolor{c2} 113 & 56         & 53             & 4, 4   & 114 & 107 \\
\hline
\rowcolor{c2} 121 & 60         & 57             & 4, 4   & 122 & 115 \\
\hline
\rowcolor{c3} 125 & 62         & 58             & 4, 4, 4   & 126 & 119 \\
\hline
\end{tabular}

\vspace{5mm}
\caption{Bounds on $\RD(\PSU(2,q))$.}
\end{table}

\newpage
\section{Explicit bounds on \texorpdfstring{$\RD(\PSU(3,q))$}{RD(PSU(3,q))}}\label{appendix:bounds_on_psu3}

We carry out the same process here on $\PSU(3,q)$ for small prime powers of $q$, using its smallest irreducible linear representation, which always lifts to a linear representation of the Schur cover $\SU(3,q)$. Take note that, as in the proof of Theorem \ref{thm:bounds_of_rd_psu3}, the algorithm derived from Theorem \ref{thm:the_game} is dominated by quartic invariants except for finitely many cases ($q=2, 3, 7, 9,$ and $13$).

\begin{table}[h!]\label{table:PSU3q_bounds}
\centering
\begin{tabular}{|c||c|c|c|c|c|}
\hline
$q$ & $\dim V$   & Bound by Thm. \ref{thm:the_game} & Degree of Invariants & $\mu(G)$ & Bound by $\RD(\mu(G))$ \\
\hline
\hline
\rowcolor{c4} 2 & 2 & 1 & None & 9 & 4  \\
\hline
\rowcolor{c3} 3 & 6 & 4 & 6 & 28 &   22\\
\hline
\rowcolor{c2} 4 & 12 & 10 & 4 & 65 &   59\\
\hline
\rowcolor{c3} 5 & 20 & 17 & 4, 4 & 50 &  44 \\
\hline
\rowcolor{c2} 7 & 42 & 39 & 4, 6 & 344 &   336\\
\hline
\rowcolor{c2} 8 & 56 & 53 & 4, 4 & 513 &  505 \\
\hline
\rowcolor{c2} 9 & 72 & 69 & 4, 6 & 730 & 722  \\
\hline
\rowcolor{c2} 11 & 110 & 106 & 4, 4, 4 & 1332 & 1324  \\
\hline
\rowcolor{c2} 13 & 156 & 152 & 4, 4, 6 & 2198 &  2189 \\
\hline
\rowcolor{c2} 16 & 240 & 236 & 4, 4, 4 & 4097 &  4088 \\
\hline
\rowcolor{c2} 17 & 272 & 267 & 4, 4, 4, 4 & 4914 &   4905\\
\hline
\rowcolor{c2} 19 & 342 & 338 & 4, 4, 4 & 6860 &   6851\\
\hline
\rowcolor{c2} 23 & 506 & 501 & 4, 4, 4, 4 & 12168 &   12159\\
\hline
\rowcolor{c2} 25 & 600 & 595 & 4, 4, 4, 4 & 15626 &   15616\\
\hline
\rowcolor{c2} 27 & 702 & 697 & 4, 4, 4, 4 & 19684 & 19674  \\
\hline
\rowcolor{c2} 29 & 812 & 807 & 4, 4, 4, 4 & 24390 & 24380  \\
\hline
\rowcolor{c2} 31 & 930 & 925 & 4, 4, 4, 4 & 29792 &  29782 \\
\hline
\rowcolor{c2} 32 & 992 & 987 & 4, 4, 4, 4 & 32769 & 32759  \\
\hline
\rowcolor{c2} 37 & 1332 & 1326 & 4, 4, 4, 4, 4 & 50654 & 50644  \\
\hline
\rowcolor{c2} 41 & 1640 & 1634 & 4, 4, 4, 4, 4 & 68922 & 68912  \\
\hline
\rowcolor{c2} 43 & 1806 & 1800 & 4, 4, 4, 4, 4 & 79508 & 79498  \\
\hline
\rowcolor{c2} 47 & 2162 & 2156 & 4, 4, 4, 4, 4 & 103824 & 103814  \\
\hline
\rowcolor{c2} 49 & 2352 & 2346 & 4, 4, 4, 4, 4 & 117650 & 117640  \\
\hline
\rowcolor{c2} 53 & 2756 & 2750 & 4, 4, 4, 4, 4 & 148878 & 148868  \\
\hline
\rowcolor{c2} 59 & 3422 & 3416 & 4, 4, 4, 4, 4 & 205380 & 205369  \\
\hline
\rowcolor{c2} 61 & 3660 & 3654 & 4, 4, 4, 4, 4 & 226982 & 226971  \\
\hline
\rowcolor{c2} 64 & 4032 & 4026 & 4, 4, 4, 4, 4 & 262145 & 262134  \\
\hline
\rowcolor{c2} 67 & 4422 & 4415 & 4, 4, 4, 4, 4, 4 & 300764 & 300753  \\
\hline
\rowcolor{c2} 71 & 4970 & 4963 & 4, 4, 4, 4, 4, 4 & 357912 & 357901  \\
\hline
\rowcolor{c2} 73 & 5256 & 5249 & 4, 4, 4, 4, 4, 4 & 389018 & 389007  \\
\hline
\rowcolor{c2} 79 & 6162 & 6155 & 4, 4, 4, 4, 4, 4 & 493040 & 493029  \\
\hline
\rowcolor{c2} 81 & 6480 & 6473 & 4, 4, 4, 4, 4, 4 & 531442 & 531431  \\
\hline
\rowcolor{c2} 83 & 6806 & 6799 & 4, 4, 4, 4, 4, 4 & 571788 & 571777  \\
\hline
\rowcolor{c2} 89 & 7832 & 7825 & 4, 4, 4, 4, 4, 4 & 704970 & 704959  \\
\hline
\rowcolor{c2} 97 & 9312 & 9305 & 4, 4, 4, 4, 4, 4 & 912674 & 912663  \\
\hline
\rowcolor{c2} 101 & 10100 & 10093 & 4, 4, 4, 4, 4, 4 & 1030302 & 1030291  \\
\hline
\rowcolor{c2} 103 & 10506 & 10499 & 4, 4, 4, 4, 4, 4 & 1092728 & 1092717  \\
\hline
\rowcolor{c2} 107 & 11342 & 11335 & 4, 4, 4, 4, 4, 4 & 1225044 & 1225033  \\
\hline
\rowcolor{c2} 109 & 11772 & 11765 & 4, 4, 4, 4, 4, 4 & 1295030 & 1295019  \\
\hline
\rowcolor{c2} 113 & 12656 & 12649 & 4, 4, 4, 4, 4, 4 & 1442898 & 1442887  \\
\hline
\rowcolor{c2} 121 & 14520 & 14513 & 4, 4, 4, 4, 4, 4 & 1771562 & 1771550   \\
\hline
\rowcolor{c2} 125 & 15500 & 15493 & 4, 4, 4, 4, 4, 4 & 1953126 & 1953114 \\
\hline
\end{tabular}

\vspace{5mm}
\caption{Bounds on $\RD(\PSU(3,q))$.}
\end{table}

\section{Powers in \texorpdfstring{$\PSU(3,q)$}{PSU(3,q)}}\label{appendix:powers_in_psu3}

Here we include the power maps needed to recover the averages of $\chi_{\Sym^k V}$ when $k = 2, 3, 4$ for the representation $V$ of $\PSU(3,q)$ discussed in Section \ref{subsec:invariant_forms}, which is used to prove Theorem \ref{thm:small_invariants}. These are computed using the class representatives given in \cite{FS1973}, where we remind the reader of corrections mentioned in \cite{Orevkov2013}. The tables are written so that the columns describe the distribution of $g \mapsto g^k$ of the indicated conjugacy class type, so each column sums to the total number of classes of that type.

\subsection{Second powers}
\;

\begin{table}[h!]
\centering
\renewcommand{\arraystretch}{1.75}
\begin{tabular}{r|c|c|c|c|c|c|c|c}
\# of Classes & $1$ & $1$ & $1$ & $q$ & $q$ & $\tfrac{q^2-q}{6}$ & $\tfrac{q^2-q}{2} - 1$ & $\tfrac{q^2-q}{3}$ \\ \hline
Type &  

$C_1$  & $C_2$  & $C_3^{\ell}$  & $C_4^k$  & $C_5^k$  & $C_6^{k,\ell,m}$  & $C_7^k$  & $C_8^k$  \\ \hline \hline 
$C_1$  &  $1$  &  $0$  &  $0$  &  $1$  &  $0$  &  $0$  &  $0$  &  $0$  \\ \hline
$C_2$  &  $0$  &  $1$  &  $0$  &  $0$  &  $1$  &  $0$  &  $0$  &  $0$  \\ \hline
$C_3^{\ell}$  &  $0$  &  $0$  &  $1$  &  $0$  &  $0$  &  $0$  &  $0$  &  $0$  \\ \hline
$C_4^k$  &  $0$  &  $0$  &  $0$  &  $q - 1$  &  $0$  &  $\frac{q}{2} - \frac{1}{2}$  &  $\frac{q}{2} + \frac{1}{2}$  &  $0$  \\ \hline
$C_5^k$  &  $0$  &  $0$  &  $0$  &  $0$  &  $q - 1$  &  $0$  &  $0$  &  $0$  \\ \hline
$C_6^{k,\ell,m}$  &  $0$  &  $0$  &  $0$  &  $0$  &  $0$  & $\frac{q^2}{6} - \frac{2q}{3} + \frac{1}{2}$  &  $0$  &  $0$  \\ \hline
$C_7^k$  &  $0$  &  $0$  &  $0$  &  $0$  &  $0$  &  $0$  & $\frac{q^2}{2} - q - \frac{3}{2}$  &  $0$  \\ \hline
$C_8^k$ \ &  $0$ \ &  $0$ \ &  $0$ \ &  $0$ \ &  $0$ \ & $0$ \ &  $0$ \ &  $\frac{q^2}{3} - \frac{q}{3}$ \\

\end{tabular}
\caption{Distribution of conjugacy classes under $g \mapsto g^2$ when $q \equiv 1,3 \bmod 6$}
\end{table}

\bigskip

\begin{table}[h!]
\centering
\renewcommand{\arraystretch}{1.75}
\begin{tabular}{r|c|c|c|c|c|c|c|c|c}
\# of Classes & $1$ & $1$ & $3$ & $\frac{q-2}{3}$ & $\frac{q-2}{3}$ & $1$ & $\tfrac{q^2-q-2}{18}$ & $\tfrac{q^2-q-2}{6}$ & $\tfrac{q^2-q-2}{9}$ \\ \hline

 Type & $C_1$  & $C_2$  & $C_3^{\ell}$  & $C_4^k$  & $C_5^k$  & $C_6'$  & $C_6^{k,\ell,m}$  & $C_7^k$  & $C_8^k$  \\ \hline \hline 
 $C_1$  &  $1$  &  $1$  &  $0$  &  $0$  &  $0$  &  $0$  &  $0$  &  $0$  &  $0$  \\ \hline 
 $C_2$  &  $0$  &  $0$  &  $3$  &  $0$  &  $0$  &  $0$  &  $0$  &  $0$  &  $0$  \\ \hline
 $C_3^{\ell}$  &  $0$  &  $0$  &  $0$  &  $0$  &  $0$  &  $0$  &  $0$  &  $0$  &  $0$  \\ \hline
 $C_4^k$  &  $0$  &  $0$  &  $0$  &  $\frac{q}{3} - \frac{2}{3}$  &  $\frac{q}{3} - \frac{2}{3}$  &  $0$  &  $0$  &  $0$  &  $0$  \\ \hline
 $C_5^k$  &  $0$  &  $0$  &  $0$  &  $0$  &  $0$  &  $0$  &  $0$  &  $0$  &  $0$  \\ \hline
 $C_6'$  &  $0$  &  $0$  &  $0$  &  $0$  &  $0$  &  $1$  &  $0$  &  $0$  &  $0$  \\ \hline
 $C_6^{k,\ell,m}$  &  $0$  &  $0$  &  $0$  &  $0$  &  $0$  &  $0$  &  $\frac{q^2}{18} - \frac{q}{18} - \frac{1}{9}$  &  $0$  &  $0$  \\ \hline
 $C_7^k$  &  $0$  &  $0$  &  $0$  &  $0$  &  $0$  &  $0$  &  $0$  &  $\frac{q^2}{6} - \frac{q}{6} - \frac{1}{3}$  &  $0$  \\ \hline
 $C_8^k$ \ &  $0$ \ &  $0$ \ &  $0$ \ &  $0$ \ &  $0$ \ &  $0$ \ &  $0$ \ &  $0$ \ &  $\frac{q^2}{9} - \frac{q}{9} - \frac{2}{9}$

\end{tabular}
\caption{Distribution of conjugacy classes under $g \mapsto g^2$ when $q \equiv 2 \bmod 6$}
\end{table}

\clearpage

\begin{table}[h!]
\centering
\renewcommand{\arraystretch}{1.75}
\begin{tabular}{r|c|c|c|c|c|c|c|c}
\# of Classes & $1$ & $1$ & $1$ & $q$ & $q$ & $\tfrac{q^2-q}{6}$ & $\tfrac{q^2-q}{2} - 1$ & $\tfrac{q^2-q}{3}$ \\
\hline
Type & $C_1$ & $C_2$ & $C_3^{\ell}$ & $C_4^k$ & $C_5^k$ & $C_6^{k,\ell,m}$ & $C_7^k$ & $C_8^k$ \\ \hline\hline
$C_1$ & $1$ & $1$ & $0$ & $0$ & $0$ & $0$ & $0$ & $0$ \\ \hline
$C_2$ & $0$ & $0$ & $1$ & $0$ & $0$ & $0$ & $0$ & $0$ \\ \hline
$C_3^{\ell}$ & $0$ & $0$ & $0$ & $0$ & $0$ & $0$ & $0$ & $0$ \\ \hline
$C_4^k$ & $0$ & $0$ & $0$ & $q$ & $q$ & $0$ & $0$ & $0$ \\ \hline
$C_5^k$ & $0$ & $0$ & $0$ & $0$ & $0$ & $0$ & $0$ & $0$ \\ \hline
$C_6^{k,\ell,m}$ & $0$ & $0$ & $0$ & $0$ & $0$ & $\tfrac{q^2}{6} - \tfrac{q}{6}$ & $0$ & $0$ \\ \hline
$C_7^k$ & $0$ & $0$ & $0$ & $0$ & $0$ & $0$ & $\tfrac{q^2}{6} -\tfrac{q}{2} - 1$ & $0$ \\ \hline
$C_8^k$ & $0$ & $0$ & $0$ & $0$ & $0$ & $0$ & $0$ & $\tfrac{q^2}{3}-\tfrac{q}{3}$
\end{tabular}
\caption{Distribution of conjugacy classes under $g \mapsto g^2$ when $q \equiv 4 \bmod 6$.}
\end{table}

\bigskip

\begin{table}[h!]
\centering
\renewcommand{\arraystretch}{1.75}
\begin{tabular}{r|c|c|c|c|c|c|c|c|c}
\# of Classes & $1$ & $1$ & $3$ & $\frac{q-2}{3}$ & $\frac{q-2}{3}$ & $1$ & $\tfrac{q^2-q-2}{18}$ & $\tfrac{q^2-q-2}{6}$ & $\tfrac{q^2-q-2}{9}$ \\ \hline
 Type & $C_1$  & $C_2$  & $C_3^{\ell}$  & $C_4^k$  & $C_5^k$  & $C_6'$  & $C_6^{k,\ell,m}$  & $C_7^k$  & $C_8^k$  \\ \hline \hline 
 $C_1$  &  $1$  &  $0$  &  $0$  &  $1$  &  $0$  &  $0$  &  $0$  &  $0$  &  $0$  \\ \hline
 $C_2$  &  $0$  &  $1$  &  $0$  &  $0$  &  $1$  &  $0$  &  $0$  &  $0$  &  $0$  \\ \hline
 $C_3^{\ell}$  &  $0$  &  $0$  &  $3$  &  $0$  &  $0$  &  $0$  &  $0$  &  $0$  &  $0$  \\ \hline 
 $C_4^k$  &  $0$  &  $0$  &  $0$  &  $\frac{q}{3} - \frac{5}{3}$  &  $0$  &  $0$  &  $\frac{q}{6} - \frac{5}{6}$  &  $\frac{q}{6} + \frac{1}{6}$  &  $0$  \\ \hline
 $C_5^k$  &  $0$  &  $0$  &  $0$  &  $0$  &  $\frac{q}{3} - \frac{5}{3}$  &  $0$  &  $0$  &  $0$  &  $0$  \\ \hline
 $C_6'$  &  $0$  &  $0$  &  $0$  &  $0$  &  $0$  &  $1$  &  $1$  &  $0$  &  $0$  \\ \hline
 $C_6^{k,\ell,m}$  &  $0$  &  $0$  &  $0$  &  $0$  &  $0$  &  $0$  &  $\frac{q^2}{18} - \frac{2q}{9} - \frac{5}{18}$  &  $0$  &  $0$  \\ \hline
 $C_7^k$  &  $0$  &  $0$  &  $0$  &  $0$  &  $0$  &  $0$  &  $0$  &  $\frac{q^2}{6} - \frac{q}{3} - \frac{1}{2}$  &  $0$  \\ \hline
 $C_8^k$ \ &  $0$ \ &  $0$ \ &  $0$ \ &  $0$ \ &  $0$ \ &  $0$ \ &  $0$ \ &  $0$ \ &  $\frac{q^2}{9} - \frac{q}{9} - \frac{2}{9}$

\end{tabular}
\caption{Distribution of conjugacy classes under $g \mapsto g^2$ when $q \equiv 5\bmod 6$}
\end{table}

\clearpage

\subsection{Third powers}
\;

\begin{table}[h!]
\centering
\renewcommand{\arraystretch}{1.75}
\begin{tabular}{r|c|c|c|c|c|c|c|c}
\# of Classes & $1$ & $1$ & $1$ & $q$ & $q$ & $\tfrac{q^2-q}{6}$ & $\tfrac{q^2-q}{2} - 1$ & $\tfrac{q^2-q}{3}$ \\ \hline
 Type & $C_1$  & $C_2$  & $C_3^{\ell}$  & $C_4^k$  & $C_5^k$  & $C_6^{k,\ell,m}$  & $C_7^k$  & $C_8^k$  \\ \hline \hline 
 $C_1$  &  $1$  &  $1$  &  $1$  &  $0$  &  $0$  &  $0$  &  $0$  &  $0$  \\ \hline 
 $C_2$  &  $0$  &  $0$  &  $0$  &  $0$  &  $0$  &  $0$  &  $0$  &  $0$  \\ \hline
 $C_3^{\ell}$  &  $0$  &  $0$  &  $0$  &  $0$  &  $0$  &  $0$  &  $0$  &  $0$  \\ \hline
 $C_4^k$  &  $0$  &  $0$  &  $0$  &  $q$  &  $q$  &  $0$  &  $0$  &  $0$  \\ \hline
 $C_5^k$  &  $0$  &  $0$  &  $0$  &  $0$  &  $0$  &  $0$  &  $0$  &  $0$  \\ \hline
 $C_6^{k,\ell,m}$  &  $0$  &  $0$  &  $0$  &  $0$  &  $0$  &  $\frac{q^2}{6} - \frac{q}{6}$  &  $0$  &  $0$  \\ \hline 
 $C_7^k$  &  $0$  &  $0$  &  $0$  &  $0$  &  $0$  &  $0$ &  $\frac{q^2}{6} - \frac{q}{2} - 1$  &  $0$  \\ \hline
 $C_8^k$ \ &  $0$ \ &  $0$ \ &  $0$ \ &  $0$ \ &  $0$ \ &  $0$ \ &  $0$ \ &  $\frac{q^2}{3} - \frac{q}{3}$ 

\end{tabular}
\caption{Distribution of conjugacy classes under $g \mapsto g^3$ when $q \equiv 0 \bmod 3$}
\end{table}

\bigskip

\begin{table}[h!]
\centering
\renewcommand{\arraystretch}{1.75}
\begin{tabular}{r|c|c|c|c|c|c|c|c}
\# of Classes & $1$ & $1$ & $1$ & $q$ & $q$ & $\tfrac{q^2-q}{6}$ & $\tfrac{q^2-q}{2} - 1$ & $\tfrac{q^2-q}{3}$ \\ \hline
 Type & $C_1$  & $C_2$  & $C_3^{\ell}$  & $C_4^k$  & $C_5^k$  & $C_6^{k,\ell,m}$  & $C_7^k$  & $C_8^k$  \\ \hline \hline 
 $C_1$  &  $1$  &  $0$  &  $0$  &  $0$  &  $0$  &  $0$  &  $1$  &  $0$  \\ \hline
 $C_2$  &  $0$  &  $1$  &  $0$  &  $0$  &  $0$  &  $0$  &  $0$  &  $0$  \\ \hline
 $C_3^{\ell}$  &  $0$  &  $0$  &  $1$  &  $0$  &  $0$  &  $0$  &  $0$  &  $0$  \\ \hline
 $C_4^k$  &  $0$  &  $0$  &  $0$  &  $q$  &  $0$  &  $0$  &  $q$  &  $0$  \\ \hline
 $C_5^k$  &  $0$  &  $0$  &  $0$  &  $0$  &  $q$  &  $0$  &  $0$  &  $0$  \\ \hline 
 $C_6^{k,\ell,m}$  &  $0$  &  $0$  &  $0$  &  $0$  &  $0$  &  $\frac{q^2}{6} - \frac{q}{6}$  &  $0$  &  $0$  \\ \hline
 $C_7^k$  &  $0$  &  $0$  &  $0$  &  $0$  &  $0$  &  $0$  &  $\frac{q^2}{2} - \frac{3q}{2} - 2$  &  $0$  \\ \hline
 $C_8^k$ \ &  $0$ \ &  $0$ \ &  $0$ \ &  $0$ \ &  $0$ \ &  $0$ \ &  $0$ \ &  $\frac{q^2}{3} - \frac{q}{3}$  

\end{tabular}
\caption{Distribution of conjugacy classes under $g \mapsto g^3$ when $q \equiv 1 \bmod 3$}
\end{table}

\clearpage

\begin{table}[h!]
\centering
\renewcommand{\arraystretch}{1.75}
\begin{tabular}{r|c|c|c|c|c|c|c|c|c}
\# of Classes & $1$ & $1$ & $3$ & $\frac{q-2}{3}$ & $\frac{q-2}{3}$ & $1$ & $\tfrac{q^2-q-2}{18}$ & $\tfrac{q^2-q-2}{6}$ & $\tfrac{q^2-q-2}{9}$ \\ \hline
 Type  & $C_1$  & $C_2$  & $C_3^{\ell}$  & $C_4^k$  & $C_5^k$  & $C_6'$  & $C_6^{k,\ell,m}$  & $C_7^k$  & $C_8^k$  \\ \hline \hline 
 $C_1$  &  $1$  &  $0$  &  $0$  &  $0$  &  $0$  &  $1$  &  $0$  &  $0$  &  $0$  \\ \hline
 $C_2$  &  $0$  &  $1$  &  $0$  &  $0$  &  $0$  &  $0$  &  $0$  &  $0$  &  $0$  \\ \hline
 $C_3^{\ell}$  &  $0$  &  $0$  &  $3$  &  $0$  &  $0$  &  $0$  &  $0$  &  $0$  &  $0$  \\ \hline
 $C_4^k$  &  $0$  &  $0$  &  $0$  &  $\frac{q}{3} - \frac{2}{3}$  &  $0$  &  $0$  &  $\frac{q}{3} - \frac{2}{3}$  &  $0$  &  $0$  \\ \hline
 $C_5^k$  &  $0$  &  $0$  &  $0$  &  $0$  &  $\frac{q}{3} - \frac{2}{3}$  &  $0$  &  $0$  &  $0$  &  $0$  \\ \hline
 $C_6'$  &  $0$  &  $0$  &  $0$  &  $0$  &  $0$  &  $0$  &  $0$  &  $0$  &  $0$  \\ \hline
 $C_6^{k,\ell,m}$  &  $0$  &  $0$  &  $0$  &  $0$  &  $0$  &  $0$  &  $\frac{q^2}{18} - \frac{7q}{18} + \frac{5}{9}$  &  $0$  &  $0$  \\ \hline
 $C_7^k$  &  $0$  &  $0$  &  $0$  &  $0$  &  $0$  &  $0$  &  $0$  &  $\frac{q^2}{6} - \frac{q}{6} - \frac{1}{3}$  &  $0$  \\ \hline
 $C_8^k$ \ &  $0$ \ &  $0$ \ &  $0$ \ &  $0$ \ &  $0$ \ &  $0$ \ &  $0$ \ &  $0$ \ &  $\frac{q^2}{9} - \frac{q}{9} - \frac{2}{9}$

\end{tabular}
\caption{Distribution of conjugacy classes under $g \mapsto g^3$ when $q \equiv 2,5 \bmod 9$}
\end{table}

\bigskip

\begin{table}[h!]
\centering
\renewcommand{\arraystretch}{1.75}
\begin{tabular}{r|c|c|c|c|c|c|c|c|c}
\# of Classes & $1$ & $1$ & $3$ & $\frac{q-2}{3}$ & $\frac{q-2}{3}$ & $1$ & $\tfrac{q^2-q-2}{18}$ & $\tfrac{q^2-q-2}{6}$ & $\tfrac{q^2-q-2}{9}$ \\ \hline
Type & $C_1$  & $C_2$  & $C_3^{\ell}$  & $C_4^k$  & $C_5^k$  & $C_6'$  & $C_6^{k,\ell,m}$  & $C_7^k$  & $C_8^k$  \\ \hline \hline 
$C_1$  &  $1$  &  $0$  &  $0$  &  $2$  &  $0$  &  $1$  &  $0$  &  $0$  &  $0$  \\ \hline
$C_2$  &  $0$  &  $1$  &  $0$  &  $0$  &  $2$  &  $0$  &  $0$  &  $0$  &  $0$  \\ \hline
$C_3^{\ell}$  &  $0$  &  $0$  &  $3$  &  $0$  &  $0$  &  $0$  &  $0$  &  $0$  &  $0$  \\ \hline
$C_4^k$  &  $0$  &  $0$  &  $0$  &  $\frac{q}{3} - \frac{8}{3}$  &  $0$  &  $0$  &  $\frac{q}{3} - \frac{8}{3}$  &  $0$  &  $0$  \\ \hline 
$C_5^k$  &  $0$  &  $0$  &  $0$  &  $0$  &  $\frac{q}{3} - \frac{8}{3}$  &  $0$  &  $0$  &  $0$  &  $0$  \\ \hline
$C_6'$  &  $0$  &  $0$  &  $0$  &  $0$  &  $0$  &  $0$  &  $3$  &  $0$  &  $0$  \\ \hline
$C_6^{k,\ell,m}$  &  $0$  &  $0$  &  $0$  &  $0$  &  $0$  &  $0$  &  $\frac{q^2}{18} - \frac{7q}{18} - \frac{4}{9}$  &  $0$  &  $0$  \\ \hline
$C_7^k$  &  $0$  &  $0$  &  $0$  &  $0$  &  $0$  &  $0$  &  $0$  &  $\frac{q^2}{6} - \frac{q}{6} - \frac{1}{3}$  &  $0$  \\ \hline
$C_8^k$ \ &  $0$ \ &  $0$ \ &  $0$ \ &  $0$ \ &  $0$ \ &  $0$ \ &  $0$ \ &  $0$ \ &  $\frac{q^2}{9} - \frac{q}{9} - \frac{2}{9}$

\end{tabular}
\caption{Distribution of conjugacy classes under $g \mapsto g^3$ when $q \equiv 8 \bmod 9$}
\end{table}

\clearpage
\subsection{Fourth powers} \;

\begin{table}[h!]
\centering
\renewcommand{\arraystretch}{1.75}
\begin{tabular}{r|c|c|c|c|c|c|c|c}
\# of Classes & $1$ & $1$ & $1$ & $q$ & $q$ & $\tfrac{q^2-q}{6}$ & $\tfrac{q^2-q}{2} - 1$ & $\tfrac{q^2-q}{3}$ \\ \hline
Type  & $C_1$  & $C_2$  & $C_3^{\ell}$  & $C_4^k$  & $C_5^k$  & $C_6^{k,\ell,m}$  & $C_7^k$  & $C_8^k$  \\ \hline  \hline 
$C_1$  &  $1$  &  $0$  &  $0$  &  $1$  &  $0$  &  $0$  &  $1$  &  $0$  \\ \hline
$C_2$  &  $0$  &  $1$  &  $0$  &  $0$  &  $1$  &  $0$  &  $0$  &  $0$  \\ \hline 
$C_3^{\ell}$  &  $0$  &  $0$  &  $1$  &  $0$  &  $0$  &  $0$  &  $0$  &  $0$  
\\ \hline $C_4^k$  &  $0$  &  $0$  &  $0$  &  $q - 1$  &  $0$  &  $\frac{q}{2} - \frac{1}{2}$  &  $\frac{3q}{2} + \frac{1}{2}$  &  $0$  
\\ \hline $C_5^k$  &  $0$  &  $0$  &  $0$  &  $0$  &  $q - 1$  &  $0$  &  $0$  &  $0$  
\\ \hline $C_6^{k,\ell,m}$  &  $0$  &  $0$  &  $0$  &  $0$  &  $0$  &  $\frac{q^2}{6} - \frac{2q}{3} + \frac{1}{2}$  &  $0$  &  $0$  
\\ \hline $C_7^k$  &  $0$  &  $0$  &  $0$  &  $0$  &  $0$  &  $0$  &  $\frac{q^2}{2} - 2q - \frac{5}{2}$  &  $0$  
\\ \hline $C_8^k$ \ &  $0$ \ &  $0$ \ &  $0$ \ &  $0$ \ &  $0$ \ &  $0$ \ &  $0$ \ &  $\frac{q^2}{3} - \frac{q}{3}$

\end{tabular}
\caption{Distribution of conjugacy classes under $g \mapsto g^4$ when $q \equiv 1, 9 \bmod 12$}
\end{table}

\bigskip

\begin{table}[h!]
\centering
\renewcommand{\arraystretch}{1.75}
\begin{tabular}{r|c|c|c|c|c|c|c|c|c}
\# of Classes & $1$ & $1$ & $3$ & $\frac{q-2}{3}$ & $\frac{q-2}{3}$ & $1$ & $\tfrac{q^2-q-2}{18}$ & $\tfrac{q^2-q-2}{6}$ & $\tfrac{q^2-q-2}{9}$ \\ \hline
 Type  & $C_1$  & $C_2$  & $C_3^{\ell}$  & $C_4^k$  & $C_5^k$  & $C_6'$  & $C_6^{k,\ell,m}$  & $C_7^k$  & $C_8^k$  \\ \hline \hline 
 $C_1$  &  $1$  &  $1$  &  $3$  &  $0$  &  $0$  &  $0$  &  $0$  &  $0$  &  $0$  \\ \hline
 $C_2$  &  $0$  &  $0$  &  $0$  &  $0$  &  $0$  &  $0$  &  $0$  &  $0$  &  $0$  \\ \hline
 $C_3^{\ell}$  &  $0$  &  $0$  &  $0$  &  $0$  &  $0$  &  $0$  &  $0$  &  $0$  &  $0$  \\ \hline
 $C_4^k$  &  $0$  &  $0$  &  $0$  &  $\frac{q}{3} - \frac{2}{3}$  &  $\frac{q}{3} - \frac{2}{3}$  &  $0$  &  $0$  &  $0$  &  $0$  \\ \hline
 $C_5^k$  &  $0$  &  $0$  &  $0$  &  $0$  &  $0$  &  $0$  &  $0$  &  $0$  &  $0$  \\ \hline
 $C_6'$  &  $0$  &  $0$  &  $0$  &  $0$  &  $0$  &  $1$  &  $0$  &  $0$  &  $0$  \\ \hline
 $C_6^{k,\ell,m}$  &  $0$  &  $0$  &  $0$  &  $0$  &  $0$  &  $0$  &  $\frac{q^2}{18} - \frac{q}{18} - \frac{1}{9}$  &  $0$  &  $0$  \\ \hline
 $C_7^k$  &  $0$  &  $0$  &  $0$  &  $0$  &  $0$  &  $0$  &  $0$  &  $\frac{q^2}{6} - \frac{q}{6} - \frac{1}{3}$  &  $0$  \\ \hline
 $C_8^k$ \ &  $0$ \ &  $0$ \ &  $0$ \ &  $0$ \ &  $0$ \ &  $0$ \ &  $0$ \ &  $0$ \ &  $\frac{q^2}{9} - \frac{q}{9} - \frac{2}{9}$

\end{tabular}
\caption{Distribution of conjugacy classes under $g \mapsto g^4$ when $q \equiv 2 \bmod 6$}
\end{table}

\clearpage

\begin{table}[h!]
\centering
\renewcommand{\arraystretch}{1.75}
\begin{tabular}{r|c|c|c|c|c|c|c|c}
\# of Classes & $1$ & $1$ & $1$ & $q$ & $q$ & $\tfrac{q^2-q}{6}$ & $\tfrac{q^2-q}{2} - 1$ & $\tfrac{q^2-q}{3}$ \\ \hline
 Type  & $C_1$  & $C_2$  & $C_3^{\ell}$  & $C_4^k$  & $C_5^k$  & $C_6^{k,\ell,m}$  & $C_7^k$  & $C_8^k$  \\ \hline \hline 
 $C_1$  &  $1$  &  $0$  &  $0$  &  $3$  &  $0$  &  $1$  &  $0$  &  $0$  \\ \hline
 $C_2$  &  $0$  &  $1$  &  $0$  &  $0$  &  $3$  &  $0$  &  $0$  &  $0$  \\ \hline
 $C_3^{\ell}$  &  $0$  &  $0$  &  $1$  &  $0$  &  $0$  &  $0$  &  $0$  &  $0$  \\ \hline
 $C_4^k$  &  $0$  &  $0$  &  $0$  &  $q - 3$  &  $0$  &  $\frac{3q}{2} - \frac{9}{2}$  &  $\frac{q}{2} + \frac{1}{2}$  &  $0$  \\ \hline
 $C_5^k$  &  $0$  &  $0$  &  $0$  &  $0$  &  $q - 3$  &  $0$  &  $0$  &  $0$  \\ \hline
 $C_6^{k,\ell,m}$  &  $0$  &  $0$  &  $0$  &  $0$  &  $0$  &  $\frac{q^2}{6} - \frac{5q}{3} + \frac{7}{2}$  &  $0$  &  $0$  \\ \hline
 $C_7^k$  &  $0$  &  $0$  &  $0$  &  $0$  &  $0$   &  $0$  &  $\frac{q^2}{2} - q - \frac{3}{2}$  &  $0$  \\ \hline
 $C_8^k$ \ &  $0$ \ &  $0$ \ &  $0$ \ &  $0$ \ &  $0$ \ &  $0$ \ &  $0$ \ &  $\frac{q^2}{3} - \frac{q}{3}$

\end{tabular}
\caption{Distribution of conjugacy classes under $g \mapsto g^4$ when $q \equiv 3, 7 \bmod 12$}
\end{table}

\bigskip

\begin{table}[h!]
\centering
\renewcommand{\arraystretch}{1.75}
\begin{tabular}{r|c|c|c|c|c|c|c|c|c}
\# of Classes & $1$ & $1$ & $1$ & $q$ & $q$ & $\tfrac{q^2-q}{6}$ & $\tfrac{q^2-q}{2} - 1$ & $\tfrac{q^2-q}{3}$ \\ \hline
 Type & $C_1$  & $C_2$  & $C_3^{\ell}$  & $C_4^k$  & $C_5^k$  & $C_6^{k,\ell,m}$  & $C_7^k$  & $C_8^k$  \\ \hline \hline 
 $C_1$  &  $1$  &  $1$  &  $1$  &  $0$  &  $0$  &  $0$  &  $0$  &  $0$  \\ \hline
 $C_2$  &  $0$  &  $0$  &  $0$  &  $0$  &  $0$  &  $0$  &  $0$  &  $0$  \\ \hline
 $C_3^{\ell}$  &  $0$  &  $0$  &  $0$  &  $0$  &  $0$  &  $0$  &  $0$  &  $0$  \\ \hline
 $C_4^k$  &  $0$  &  $0$  &  $0$  &  $q$  &  $q$  &  $0$  &  $0$  &  $0$  \\ \hline
 $C_5^k$  &  $0$  &  $0$  &  $0$  &  $0$  &  $0$ &  $0$  &  $0$  &  $0$  \\ \hline
 $C_6'$  &  $0$  &  $0$  &  $0$  &  $0$  &  $0$  &  $0$  &  $0$  &  $0$  \\ \hline
 $C_6^{k,\ell,m}$  &  $0$  &  $0$  &  $0$  &  $0$  &  $0$  &  $\frac{q^2}{6} - \frac{q}{6}$  &  $0$  &  $0$  \\ \hline
 $C_7^k$  &  $0$  &  $0$  &  $0$  &  $0$  &  $0$  &  $0$  &  $\frac{q^2}{2} - \frac{q}{2} - 1$  &  $0$  \\ \hline
 $C_8^k$ \ &  $0$ \ &  $0$ \ &  $0$ \ &  $0$ \ &  $0$ \ &  $0$ \ &  $0$ \ &  $\frac{q^2}{3} - \frac{q}{3}$

\end{tabular}
\caption{Distribution of conjugacy classes under $g \mapsto g^4$ when $q \equiv 4 \bmod 6$}
\end{table}

\clearpage

\begin{table}[h!]
\centering
\renewcommand{\arraystretch}{1.75}
\begin{tabular}{r|c|c|c|c|c|c|c|c|c}
\# of Classes & $1$ & $1$ & $3$ & $\frac{q-2}{3}$ & $\frac{q-2}{3}$ & $1$ & $\tfrac{q^2-q-2}{18}$ & $\tfrac{q^2-q-2}{6}$ & $\tfrac{q^2-q-2}{9}$ \\ \hline
Type & $C_1$  & $C_2$  & $C_3^{\ell}$  & $C_4^k$  & $C_5^k$  & $C_6'$  & $C_6^{k,\ell,m}$  & $C_7^k$  & $C_8^k$  \\ \hline \hline 
$C_1$  &  $1$  &  $0$  &  $0$  &  $1$  &  $0$  &  $0$  &  $0$  &  $1$  &  $0$  \\ \hline
$C_2$  &  $0$  &  $1$  &  $0$  &  $0$  &  $1$  &  $0$  &  $0$  &  $0$  &  $0$  \\ \hline
$C_3^{\ell}$  &  $0$  &  $0$  &  $3$  &  $0$  &  $0$  &  $0$  &  $0$  &  $0$  &  $0$  \\ \hline
$C_4^k$  &  $0$  &  $0$  &  $0$  &  $\frac{q}{3} - \frac{5}{3}$  &  $0$  &  $0$  &  $\frac{q}{6} - \frac{5}{6}$  &  $\frac{q}{2} - \frac{1}{2}$  &  $0$  \\ \hline
$C_5^k$  &  $0$  &  $0$  &  $0$  &  $0$  &  $\frac{q}{3} - \frac{5}{3}$  &  $0$  &  $0$  &  $0$  &  $0$  \\ \hline
$C_6'$  &  $0$  &  $0$  &  $0$  &  $0$  &  $0$  &  $1$  &  $1$  &  $0$  &  $0$  \\ \hline 
$C_6^{k,\ell,m}$  &  $0$  &  $0$  &  $0$  &  $0$  &  $0$  &  $0$  &  $\frac{q^2}{18} - \frac{2q}{9} - \frac{5}{18}$  &  $0$  &  $0$  \\ \hline
$C_7^k$  &  $0$  &  $0$  &  $0$  &  $0$  &  $0$  &  $0$  &  $0$  &  $\frac{q^2}{6} - \frac{2q}{3} - \frac{5}{6}$  &  $0$  \\ \hline
$C_8^k$ \ &  $0$ \ &  $0$ \ &  $0$ \ &  $0$ \ &  $0$ \ &  $0$ \ &  $0$ \ &  $0$ \ &  $\frac{q^2}{9} - \frac{q}{9} - \frac{2}{9}$

\end{tabular}
\caption{Distribution of conjugacy classes under $g \mapsto g^4$ when $q \equiv 5 \bmod 12$}
\end{table}

\bigskip

\begin{table}[h!]
\centering
\renewcommand{\arraystretch}{1.75}
\begin{tabular}{r|c|c|c|c|c|c|c|c|c}
\# of Classes & $1$ & $1$ & $3$ & $\frac{q-2}{3}$ & $\frac{q-2}{3}$ & $1$ & $\tfrac{q^2-q-2}{18}$ & $\tfrac{q^2-q-2}{6}$ & $\tfrac{q^2-q-2}{9}$ \\ \hline
Type & $C_1$  & $C_2$  & $C_3^{\ell}$  & $C_4^k$  & $C_5^k$  & $C_6'$  & $C_6^{k,\ell,m}$  & $C_7^k$  & $C_8^k$  \\ \hline \hline 
$C_1$  &  $1$  &  $0$  &  $0$  &  $3$  &  $0$  &  $0$  &  $1$  &  $0$  &  $0$  \\ \hline
$C_2$  &  $0$  &  $1$  &  $0$  &  $0$  &  $3$  &  $0$  &  $0$  &  $0$  &  $0$  \\ \hline
$C_3^{\ell}$  &  $0$  &  $0$  &  $3$  &  $0$  &  $0$  &  $0$  &  $0$  &  $0$  &  $0$  \\ \hline
$C_4^k$  &  $0$  &  $0$  &  $0$  &  $\frac{q}{3} - \frac{11}{3}$  &  $0$  &  $0$  &  $\frac{q}{2} - \frac{11}{2}$  &  $\frac{q}{6} + \frac{1}{6}$  &  $0$  \\ \hline
$C_5^k$  &  $0$  &  $0$  &  $0$  &  $0$  &  $\frac{q}{3} - \frac{11}{3}$  &  $0$  &  $0$  &  $0$  &  $0$  \\ \hline $C_6'$  &  $0$  &  $0$  &  $0$  &  $0$  &  $0$  &  $1$  &  $5$  &  $0$  &  $0$  \\ \hline 
$C_6^{k,\ell,m}$  &  $0$  &  $0$  &  $0$  &  $0$  &  $0$  &  $0$  &  $\frac{q^2}{18} - \frac{5q}{9} - \frac{11}{18}$  &  $0$  &  $0$  \\ \hline 
$C_7^k$  &  $0$  &  $0$  &  $0$  &  $0$  &  $0$  &  $0$  &  $0$  &  $\frac{q^2}{6} - \frac{q}{3} - \frac{1}{2}$  &  $0$  \\ \hline
$C_8^k$ \ &  $0$ \ &  $0$ \ &  $0$ \ &  $0$ \ &  $0$ \ &  $0$ \ &  $0$ \ &  $0$ \ &  $\frac{q^2}{9} - \frac{q}{9} - \frac{2}{9}$

\end{tabular}
\caption{Distribution of conjugacy classes under $g \mapsto g^4$ when $q \equiv 11 \bmod 12$}
\end{table}

\bigskip

\clearpage

\begin{bibsection}
\begin{biblist}

\bib{AS1976}{article}{
   author={Arnold, V.},
   author={Shimura, G.},
   title={Superpositions of algebraic functions},
   journal={Proc. Symposia in Pure Math.},
   volume={28},
   date={1976},
   pages={45--46},
}

\bib{Brauer1975}{article}{
   author={Brauer, Richard},
   title={On the resolvent problem},
   journal={Annali di Matematica Pura ed Applicata},
   volume={102},
   number={1},
   date={1975},
   pages={45--55},
}

\bib{Bring1786}{article}{
   label={Bri1786},
   author={Bring, E.},
   title={Meletemata qu{\ae}dam mathematica circa transformationem {\ae}quationem algebraicarum},
   journal={Typis Berlingianis},
   date={1786},
}

\bib{FW2019}{article}{
   author={Farb, Benson},
   author={Wolfson, Jesse},
   title={Resolvent degree, Hilbert's 13th problem and geometry},
   journal={L'Enseignement Math\'ematique},
   volume={65},
   number={3},
   date={2019},
   pages={303--376},
}

\bib{FKW2023}{article}{
   author={Farb, Benson},
   author={Kisin, Mark},
   author={Wolfson, Jesse},
   title={Modular functions and resolvent problems: With an appendix by Nate Harman},
   journal={Mathematische Annalen},
   volume={386},
   number={1},
   date={2023},
   pages={113--150},
}

\bib{FS1973}{article}{
   author={Frame, J. Sutherland},
   author={Simpson, William A.},
   title={The character tables for $\SL(3, q)$, $\SU(3,q^2)$, $\PSL(3,q)$, $\PSU(3,q^2)$},
   journal={Canadian Journal of Mathematics},
   volume={25},
   date={1973},
   pages={486--494},
}

\bib{GAP2024}{manual}{
   label={GAP24},
   author={The GAP Group},
   title={GAP -- Groups, Algorithms, and Programming, Version 4.14.0},
   date={2024},
   note={\url{https://www.gap-system.org}},
}

\bib{GGSW2024}{article}{
   author={G\'omez-Gonz\'ales, C.},
   author={Sutherland, Alexander},
   author={Wolfson, Jesse},
   title={Generalized versality, special points, and resolvent degree for the sporadic groups},
   journal={Journal of Algebra},
   volume={647},
   date={2024},
   pages={758--793},
}

\bib{GM2020}{book}{
   author={Geck, Meinolf},
   author={Malle, Gunter},
   title={The Character Theory of Finite Groups of Lie Type: A Guided Tour},
   series={London Mathematical Society Lecture Note Series},
   volume={187},
   publisher={Cambridge University Press},
   date={2020},
}

\bib{GMPS2015}{article}{
   author={Guest, Simon},
   author={Morris, Joy},
   author={Praeger, Cheryl},
   author={Spiga, Pablo},
   title={On the maximum orders of elements of finite almost simple groups and primitive permutation groups},
   journal={Transactions of the American Mathematical Society},
   volume={367},
   number={11},
   date={2015},
   pages={7665--7694},
}

\bib{Hamilton1836}{article}{
   label={Ham1836},
   author={Hamilton, William Rowan},
   title={Inquiry into the validity of a method recently proposed by George B. Jerrard, Esq., for transforming and resolving equations of elevated degrees},
   journal={Report of the Sixth Meeting of the British Association for the Advancement of Science},
   date={1836},
   pages={295--348},
}

\bib{Jordan1907}{article}{
   label={Jor1907},
   author={Jordan, Herbert E.},
   title={Group-characters of various types of linear groups},
   journal={American Journal of Mathematics},
   volume={29},
   number={4},
   date={1907},
   pages={387--405},
}

\bib{Karpilovsky1987}{book}{
   author={Karpilovsky, Gregory},
   title={The Schur Multiplier},
   series={London Mathematical Society Monographs},
   volume={2},
   publisher={Oxford University Press},
   date={1987},
}

\bib{Klein1878}{article}{
   label={Kle1878},
   author={Klein, Felix},
   title={Ueber die Transformation siebenter Ordnung der elliptischen Functionen},
   journal={Mathematische Annalen},
   volume={14},
   number={3},
   date={1878},
   pages={428--471},
}

\bib{HeberleSutherland2023}{article}{
   author={Heberle, C.},
   author={Sutherland, Alexander},
   title={Upper bounds on resolvent degree via Sylvester's obliteration algorithm},
   journal={New York Journal of Mathematics},
   volume={29},
   date={2023},
   pages={107--146},
}

\bib{Isaacs1994}{book}{
   author={Isaacs, I. Martin},
   title={Character Theory of Finite Groups},
   series={Dover Books on Mathematics},
   publisher={Courier Corporation},
   date={1994},
}

\bib{Klein1884}{book}{
   label={Kle1884},
   author={Klein, Felix},
   title={Vorlesungen \"uber das Ikosaeder und die Aufl\"osung der Gleichungen vom f\"unften Grade},
   publisher={Teubner},
   date={1884},
}

\bib{Orevkov2013}{article}{
   author={Orevkov, S. Yu.},
   title={Products of conjugacy classes in finite unitary groups $\operatorname{GU}(3, q^2)$ and $\operatorname{SU}(3, q^2)$},
   journal={Annales de la Facult\'e des sciences de Toulouse: Math\'ematiques},
   volume={22},
   number={2},
   date={2013},
   pages={219--251},
}

\bib{Reichstein2025}{article}{
   author={Reichstein, Z.},
   title={Hilbert's 13th problem for algebraic groups},
   journal={L'Enseignement Math\'ematique},
   volume={71},
   number={1},
   date={2025},
   pages={139--192},
}

\bib{Sage2025}{manual}{
   label={Sag25},
   author={The Sage Developers},
   title={SageMath. The Sage Mathematics Software System, Version 10.7},
   date={2025},
   note={\url{https://www.sagemath.org}},
}

\bib{Schur1907}{article}{
   label={Sch1907},
   author={Schur, J.},
   title={Untersuchungen \"uber die Darstellung der endlichen Gruppen durch gebrochene lineare Substitutionen},
   journal={Journal f\"ur die reine und angewandte Mathematik},
   volume={132},
   date={1907},
   pages={85--137},
}

\bib{Sutherland2021}{article}{
   author={Sutherland, Alexander J.},
   title={Upper bounds on resolvent degree and its growth rate},
   journal={arXiv preprint},
   volume={arXiv:2107.08139},
   date={2021},
}

\bib{Sutherland2023}{article}{
   author={Sutherland, Alexander J.},
   title={A Summary of Known Bounds on the Essential Dimension and Resolvent Degree of Finite Groups},
   journal={arXiv preprint},
   volume={arXiv:2312.04430},
   date={2023},
}

\bib{Tschirnhaus1683}{article}{
   label={Tsc1683},
   author={von Tschirnhaus, E.},
   title={Methodus auferendi omnes terminos intermedios ex data aequatione},
   journal={Acta Eruditorum},
   date={1683},
   pages={204--207},
}

\bib{Wolfson2020}{article}{
   author={Wolfson, Jesse},
   title={Tschirnhaus transformations after Hilbert},
   journal={L'Enseignement Math\'ematique},
   volume={66},
   number={3},
   date={2020},
   pages={489--540},
}
\end{biblist}
\end{bibsection}
\end{document}